\documentclass[a4paper,12pt,amsart,frenchb]{article}

\usepackage{amsmath,amsbsy,amsfonts,amssymb}
\usepackage[french]{babel}
\newcommand{\ass}[2]{\vskip0.3cm\noindent
{\bf {#1}}. { \sl {#2}}\vskip0.3cm\noindent
}

\begin{document}
 
   \title{  Caract\`eres automorphes d'un groupe r\'eductif}
\author{J.-L. Waldspurger}
\date{25 ao\^ut 2016}
\maketitle

\sl {\bf Abstract} Let $G$ be a reductive group defined over a number field. Denote $Z(\hat{G})$ the center of the dual group. Langlands has defined some homomorphism from some cohomology group of $Z(\hat{G})$ into the group of automorphic characters of $G$. We prove that it is bijective.\rm

 \bigskip
 Soit $F$ un corps de caract\'eristique $0$ qui est soit un corps local, soit un corps de nombres. On fixe une cl\^oture alg\'ebrique $\bar{F}$ de $F$. On note $\Gamma_{F}$ le groupe de Galois $Gal(\bar{F}/F)$ et $W_{F}$ le groupe de Weil de $F$. Si $F$ est un corps de nombres, on note ${\mathbb A}_{F}$ son anneau des ad\`eles et $Val(F)$ l'ensemble de ses valuations. Soit $G$ un groupe r\'eductif connexe d\'efini sur $F$. On introduit le groupe dual complexe $\hat{G}$ et son centre $Z(\hat{G})$. Ces groupes sont munis d'une action alg\'ebrique de $\Gamma_{F}$. 
 
 Si $F$ est local, on d\'efinit apr\`es Langlands un homomorphisme
 $$\alpha_{F}:H^1(W_{F};Z(\hat{G}))\to Hom_{cont}(G(F),{\mathbb C}^{\times}),$$
 l'indice "cont" indiquant qu'il s'agit d'homomorphismes continus.
 
 Si $F$ est un corps de nombres, on note $Hom_{cont}(G(F)\backslash G({\mathbb A}_{F}), {\mathbb C}^{\times})$ le groupe des homomorphismes continus de $G({\mathbb A}_{F})$ dans ${\mathbb C}^{\times}$ qui sont \'egaux \`a $1$ sur $G(F)$. D'autre part, on a un homomorphisme de localisation
 $$H^1(W _{F};Z(\hat{G}))\to \prod_{v\in Val(F)}H^1(W_{F_{v}};Z(\hat{G}))$$
 dont on note $ker^1(W_{F};Z(\hat{G}))$ le noyau. On d\'efinit alors un homomorphisme
 $$\alpha_{F}:H^1(W_{F};Z(\hat{G}))/ker^1(W_{F};Z(\hat{G}))\to Hom_{cont}(G(F)\backslash G({\mathbb A}_{F}),{\mathbb C}^{\times}).$$
 
 Si $F={\mathbb R}$, l'homomorphisme $\alpha_{F}$ n'est en g\'en\'eral ni surjectif, ni injectif. Si $F$ est $p$-adique, $\alpha_{F}$ est toujours injectif. J'ai \`a plusieurs reprises commis l'erreur d'affirmer que $\alpha_{F}$ \'etait toujours surjectif. J.-P. Labesse m'a signal\'e cette erreur et indiqu\'e la r\'ef\'erence \cite{LL} qui d\'ecrit pr\'ecis\'ement ce qu'il en est (disons toutefois que $\alpha_{F}$ est surjectif pour beaucoup de groupes usuels). Je l'en remercie vivement. Cependant, la r\'ef\'erence \cite{LL} traite aussi le cas des corps de nombres et, dans ce cas, son r\'esultat n'est pas optimal. En effet, on a le lemme suivant.
 
 \ass{Lemme}{Si $F$ est un corps de nombres, $\alpha_{F}$ est bijectif.}
 
 Peut-\^etre qu'une preuve de ce r\'esultat se trouve d\'ej\`a dans la litt\'erature. Faute de l'avoir trouv\'ee, je vais en donner une.
 
 {\bf Remarque.} Les auteurs de \cite{LL} ne supposent pas que $F$ est de caract\'eristique nulle. Par prudence, je me limiterai \`a ce cas.
 
 \bigskip
 
 Preuve. Soient $T$ et $T'$ deux tores d\'efinis sur $F$ et $f:T\to T'$ un homomorphisme d\'efini sur $F$. On d\'efinit divers groupes de cohomologie associ\'es au $2$-complexe $T\stackrel{f}{\to}T'$. F\^acheusement, la num\'erotation de ces groupes diff\`ere selon les auteurs. Nous notons $H^{i+1,i}$ le groupe dont les \'el\'ements sont des classes d'\'equivalence de paires de cocha\^{\i}nes, la premi\`ere \'etant de degr\'e $i+1$ et la seconde de degr\'e $i$. Ainsi, on a des groupes $H^{1,0}(F;T\stackrel{f}{\to}T')$, $H^{1,0}({\mathbb A}_{F};T\stackrel{f}{\to}T')$, $H^{1,0}({\mathbb A}_{F}/F;T\stackrel{f}{\to}T')$, cf. \cite{L} 1.4 ou \cite{KS} appendice C1. Ils s'inscrivent dans une suite exacte
 $$H^{1,0}(F;T\stackrel{f}{\to}T')\to H^{1,0}({\mathbb A}_{F};T\stackrel{f}{\to}T')\to H^{1,0}({\mathbb A}_{F}/F;T\stackrel{f}{\to}T').$$
 Le groupe $H^{1,0}({\mathbb A}_{F};T\stackrel{f}{\to}T')$ est un produit restreint des groupes $H^{1,0}(F_{v};T\stackrel{f}{\to}T')$ sur les places $v\in Val(F)$. Tous ces groupes sont munis de topologies qui en font des groupes localement compacts. Les fl\`eches de la suite ci-dessus sont continues.
 Dualement, on a un complexe $\hat{T}'\stackrel{\hat{f}}{\to}\hat{T}$ et on d\'efinit le groupe de cohomologie  $H^{1,0}(W_{F}; \hat{T}'\stackrel{\hat{f}}{\to}\hat{T})$ et, pour toute place $v\in Val(F)$, le groupe $H^{1,0}(W_{F_{v}};\hat{T}'\stackrel{\hat{f}}{\to}\hat{T})$. On a un homomorphisme
 $$H^{1,0}(W_{F}; \hat{T}'\stackrel{\hat{f}}{\to}\hat{T})\to \prod_{v\in Val(F)}H^{1,0}(W_{F_{v}};\hat{T}'\stackrel{\hat{f}}{\to}\hat{T}).$$
 On dispose d'un accouplement
 $$H^{1,0}(W_{F}; \hat{T}'\stackrel{\hat{f}}{\to} \hat{T})\times H^{1,0}({\mathbb A}_{F}/F;T\stackrel{f}{\to}T')\to {\mathbb C}^{\times}$$
 et, pour tout $v\in Val(F)$, d'un accouplement
 $$H^{1,0}(W_{F_{v}};\hat{T}'\stackrel{\hat{f}}{\to}\hat{T})\to H^{1,0}(F_{v};T\stackrel{f}{\to}T').$$
 Ils sont compatibles avec les homomorphismes d\'efinis ci-dessus.
 
 Notons $G_{SC}$ le rev\^etement simplement connexe du groupe d\'eriv\'e de $G$ et notons $\pi:G_{SC}\to G$ l'homomorphisme naturel. Choisissons un sous-tore maximal $T$ de $G$ d\'efini sur $F$, notons $T_{sc}$ son image r\'eciproque par $\pi$. Labesse d\'efinit les groupes $H^0_{ab}(F;G)$, $H^0_{ab}({\mathbb A}_{F};G)$ et $H^0_{ab}({\mathbb A}_{F}/F;G)$ comme \'etant $H^0_{ab}(F;T_{sc}\stackrel{\pi}{\to}T)$, etc... Ils ne d\'ependent pas du choix de $T$ car on voit que l'on peut remplacer dans leurs d\'efinitions le complexe de tores $T_{sc}\stackrel{\pi}{\to}T$ par le complexe de groupes diagonalisables $Z(G_{SC})\stackrel{\pi}{\to}Z(G)$ (o\`u $Z(G)$ est le centre de $G$ et $Z(G_{SC})$ celui de $G_{SC}$). On dispose d'homomorphismes $ab_{F}:G(F)\to H^0_{ab}(F;G)$ et $ab_{F_{v}}:G(F_{v})\to H^0(F_{v};G)$ pour tout $v\in Val(F)$, ces derniers se regroupant en un homomorphisme $ab_{{\mathbb A}_{F}}:G({\mathbb A}_{F})\to H^0_{ab}({\mathbb A}_{F};G)$. Rappelons leur d\'efinition, par exemple dans le cas de $ab_{F}$. Pour $g\in G(F)$, on choisit $g_{sc}\in G_{SC}(\bar{F})$ et $z\in Z(G)(\bar{F})$ tels que $g=\pi(g_{sc})z$. On d\'efinit un cocycle $\mu:\Gamma_{F}\to G_{SC}$ par $\mu(\sigma)=g_{sc}\sigma(g_{sc})^{-1}$ pour $\sigma\in \Gamma_{F}$. On voit qu'il prend ses valeurs dans $Z(G_{SC})(\bar{F})$ et que le couple de cocha\^{\i}nes $(\mu,z)$ est un cocycle dont la classe dans $H^0_{ab}(F;G)$ est, par d\'efinition, $ab_{F}(g)$. On voit que le noyau de $ab_{F}$ est \'egal \`a $\pi(G_{SC}(F))$ et que, pour toute place $v$, le noyau de $ab_{F_{v}}$ est $\pi(G_{SC}(F_{v}))$. On voit aussi que l'homomorphisme $ab_{{\mathbb A}_{F}}$ est continu et ouvert.

 Au tore $T$ est associ\'e un sous-tore maximal $\hat{T}$ de $\hat{G}$. Le groupe dual de $G_{SC}$ est le groupe adjoint $\hat{G}_{AD}=\hat{G}/Z(\hat{G})$. Notons $\hat{\pi}:\hat{G}\to \hat{G}_{AD}$ l'homomorphisme naturel et $\hat{T}_{ad}$ l'image de $\hat{T}$ par $\pi_{ad}$. Le complexe dual de $T_{sc}\stackrel{\pi}{\to}T$ est alors $\hat{T}\stackrel{\hat{\pi}}{\to}\hat{T}_{ad}$. On voit que les groupes de cohomologie associ\'es ce complexe sont \'egaux \`a ceux associ\'es au complexe $Z(\hat{G})\to \{1\}$, autrement dit ce sont les groupes $H^1(W_{F};Z(\hat{G}))$, etc...
 
   Nous avons ainsi introduit les objets et les homomorphismes du diagramme suivant:
   $$(1)\qquad \begin{array}{ccccc} G_{SC}(F)&&G_{SC}({\mathbb A}_{F})&&\\ \,\,\downarrow \pi&&\,\,\downarrow \pi&&\\   
   G(F)&\to&G({\mathbb A}_{F})&&\\ \quad \downarrow ab_{F}&&\quad \downarrow ab_{{\mathbb A}_{F}}&&\\ H^0_{ab}(F;G)&\stackrel{j}{\to}&H^0_{ab}({\mathbb A}_{F};G)&\stackrel{k}{\to}& H^0_{ab}({\mathbb A}_{F}/F;G)\\&&&&\\  && \prod_{v\in Val_{F}}H^{1}(W_{F_{v}};Z(\hat{G}))&\leftarrow&H^1(W_{F};Z(\hat{G}))\\ \end{array}$$
   Dans notre cas, les accouplements \'evoqu\'es plus haut sont des dualit\'es parfaites (\cite{KS} lemmes A.3.B et C.2.C), c'est-\`a-dire qu'il s'en d\'eduit un isomorphisme
   $$\alpha'_{F}:H^1(W_{F};Z(\hat{G}))\to Hom_{cont}(H^0_{ab}({\mathbb A}_{F}/F;G),{\mathbb C}^{\times})$$
   et, pour tout $v\in Val(F)$, un isomorphisme
   $$\alpha'_{F_{v}}:H^1(W_{F_{v}};Z(\hat{G}))\to Hom_{cont}(H^0_{ab}(F_{v};G),{\mathbb C}^{\times}).$$
   Ces isomorphismes sont compatibles avec les homomorphismes du diagramme (1). Parce que les $\alpha'_{F_{v}}$ sont des isomorphismes, l'image par $\alpha'_{F}$ de $ker^1(W_{F};Z(\hat{G}))$ est le sous-groupe des \'el\'ements de $Hom_{cont}(H^0_{ab}({\mathbb A}_{F}/F;G),{\mathbb C}^{\times})$ qui sont triviaux sur l'image de l'homomorphisme $k$ du diagramme (1). Or, d'apr\`es \cite{KS} lemme C.3.A, $k$ se quotiente en un hom\'eomorphisme de $H^0_{ab}({\mathbb A}_{F};G)/j(H^0_{ab}(F;G))$ sur un sous-groupe ouvert de $H^0_{ab}({\mathbb A}_{F}/F;G)$. Donc le groupe $Hom_{cont}(H^0_{ab}({\mathbb A}_{F}/F;G),{\mathbb C}^{\times})$, quotient\'e par le sous-groupe des \'el\'ements qui sont triviaux sur l'image de $k$, s'identifie \`a $Hom_{cont}(H^0_{ab}({\mathbb A}_{F};G)/j(H^0_{ab}(F;G)),{\mathbb C}^{\times})$. De $\alpha'_{F}$ se d\'eduit donc un isomorphisme
   $$\alpha''_{F}:H^1(W_{F};Z(\hat{G}))/ker^1(W_{F};Z(\hat{G}))\to Hom_{cont}(H^0_{ab}({\mathbb A}_{F};G)/j(H^0_{ab}(F;G)),{\mathbb C}^{\times}).$$
   D'autre part, de l'homomorphisme $ab_{{\mathbb A}_{F}}$ se d\'eduit un homomorphisme
   $$ab'_{{\mathbb A}_{F}}: G({\mathbb A}_{F})\to H^0_{ab}({\mathbb A}_{F};G)/j(H^0_{ab}(F;G)),$$
   dont le noyau contient $G(F)$. On a donc dualement un homomorphisme
   $$\beta:Hom_{cont}(H^0_{ab}({\mathbb A}_{F};G)/j(H^0_{ab}(F;G)),{\mathbb C}^{\times})\to Hom_{cont}(G(F)\backslash G({\mathbb A}_{F}),{\mathbb C}^{\times}).$$
   par d\'efinition, l'application $\alpha_{F}$ de l'\'enonc\'e est \'egale \`a $\beta\circ \alpha''_{F}$. Puisque $\alpha''_{F}$ est bijectif, on doit prouver que $\beta$ l'est aussi.
   
   Montrons d'abord que 
   
   (2) si $G$ est simplement connexe, $Hom_{cont}(G(F)\backslash G({\mathbb A}_{F}),{\mathbb C}^{\times})=\{1\}$.
 
  Notons $S$ l'ensemble des places $v\in Val(F)$ telles que $G$ ne soit pas quasi-d\'eploy\'e sur $F_{v}$. C'est un ensemble fini. D'apr\`es  \cite{LL}, $\alpha_{F_{v}}$ est un isomorphisme pour $v\in Val(F)-S$. Or $Z(\hat{G})=\{1\}$ puisque $\hat{G}$ est adjoint. Donc $Hom_{cont}(G(F_{v}),{\mathbb C}^{\times})=\{1\}$ pour $v\in Val(F)- S$. Un \'el\'ement $\omega$ de  $Hom_{cont}(G(F)\backslash G({\mathbb A}_{F}),{\mathbb C}^{\times})$ est donc trivial sur $G(F_{v})$ pour  ces places $v$, donc $\omega$ est trivial sur $G(F)G({\mathbb A}_{F}^S)$, o\`u ${\mathbb A}_{F}^S$ est le produit restreint des $F_{v}$ pour $v\in Val(F)-S$. Or, parce que $G$ est simplement connexe, cet ensemble  est dense dans $G({\mathbb A}_{F})$ (th\'eor\`eme de Kneser-Harder, cf. \cite{S} th\'eor\`eme 3.1). Donc $\omega$ est trivial, ce qui prouve (2). 
  
   Revenons \`a $G$ quelconque. Soit $\omega\in Hom_{cont}(G(F)\backslash G({\mathbb A}_{F}),{\mathbb C}^{\times})$. Alors $\omega\circ \pi$ est un \'el\'ement de $Hom_{cont}(G_{SC}(F)\backslash G_{SC}({\mathbb A}_{F}),{\mathbb C}^{\times})$, lequel est trivial d'apr\`es (2). Donc $\omega$ est trivial sur $\pi(G_{SC}({\mathbb A}_{F}))$. On obtient que l'application naturelle
   $$G({\mathbb A}_{F})\to G(F)\pi(G_{SC}({\mathbb A}_{F}))\backslash G({\mathbb A}_{F})$$
   se dualise en une bijection
   $$(3) \qquad Hom_{cont}(G(F)\pi(G_{SC}({\mathbb A}_{F}))\backslash G({\mathbb A}_{F}),{\mathbb C}^{\times})\to Hom_{cont}(G(F)\backslash G({\mathbb A}_{F}),{\mathbb C}^{\times}).$$
   L'homomorphisme $ab'_{{\mathbb A}_{F}}$ introduit plus haut se quotiente en un homomorphisme
   $$b:G(F)\pi(G_{SC}({\mathbb A}_{F}))\backslash G({\mathbb A}_{F})\to H^0_{ab}({\mathbb A}_{F};G)/j(H^0_{ab}(F;G))$$
   et, modulo la bijection (3), $\beta$ n'est autre que l'homomorphisme dual de $b$. D'autre part, parce que $ab_{{\mathbb A}_{F}}$ est continu et ouvert, $b$ l'est aussi. Pour montrer que $\beta$ est bijectif, il suffit de prouver que $b$ l'est. 
   
  Remarquons que, puisque les suites verticales du diagramme (1) sont exactes, le groupe $G(F)\pi(G_{SC}({\mathbb A}_{F}))\backslash G({\mathbb A}_{F})$ s'identifie avec $ab_{{\mathbb A}_{F}}(G({\mathbb A}_{F}))/j\circ ab_{F}(G(F))$ et $b$  devient l'application
  $$b:ab_{{\mathbb A}_{F}}(G({\mathbb A}_{F}))/j\circ ab_{F}(G(F))\to H^0_{ab}({\mathbb A}_{F};G)/j(H^0_{ab}(F;G))$$
  d\'eduite de l'injection  $ab_{{\mathbb A}_{F}}(G({\mathbb A}_{F}))\subset H^0_{ab}({\mathbb A}_{F};G)$.
  
   Montrons que
   
   (4) $b$ est surjectif.
   
   Notons $Val_{\infty}(F)$ l'ensemble des places archim\'ediennes de $F$. Pour toute place $v\in Val(F)$, l'image de $ab_{F_{v}}$ est un sous-groupe ouvert d'indice fini de $H^0_{ab}(F_{v};G)$. Si $v\not\in Val_{\infty}(F)$, $ab_{F_{v}}$ est surjectif, cf. \cite{L} proposition 1.6.7. Donc l'image de $b$ est un sous-groupe ouvert d'indice fini de $H^0_{ab}({\mathbb A}_{F};G)/j(H^0_{ab}(F;G))$ qui contient l'image naturelle de $H^0_{ab}({\mathbb A}_{F}^{Val_{\infty}(F)};G)$. D'apr\`es \cite{KS} lemme C.5.A, la projection de $j (H^0_{ab}(F;G))$ dans le produit $\prod_{v\in Val_{\infty}(F)}H^0_{ab}(F_{v};G)$ est d'image dense. Il revient au m\^eme de dire que le produit de $j(H^0_{ab}(F;G))$ et de l'image naturelle de $H^0_{ab}({\mathbb A}_{F}^{Val_{\infty}(F)};G)$ dans    $H^0_{ab}({\mathbb A}_{F};G)$ est dense dans ce dernier groupe, ou encore que l'image naturelle de $H^0_{ab}({\mathbb A}_{F}^{Val_{\infty}(F)};G)$ dans $H^0_{ab}({\mathbb A}_{F};G)/j(H^0_{ab}(F;G))$ est dense. Ainsi, l'image de $b$ est un sous-groupe ouvert dense, donc $b$ est surjectif.
   
   Montrons enfin que
   
   (5) $b$ est injectif. 

Il revient au m\^eme de prouver l'\'egalit\'e

(6) $j(H^0_{ab}(F;G))\cap ab_{{\mathbb A}_{F}}(G({\mathbb A}_{F}))= j\circ ab_{F}(G(F))$.

L'inclusion du membre de droite dans celui de gauche est claire et est d'ailleurs utilis\'ee dans les constructions ci-dessus. Montrons l'inclusion inverse.  L'homomorphisme $j$ s'ins\`ere dans un diagramme commutatif
$$\begin{array}{ccccccc} G_{SC}(F)&\stackrel{\pi}{\to}&G(F)&\stackrel{ab_{F}}{\to} &H^0_{ab}(F;G)&\to& H^1(F;G_{SC})\\ \downarrow&&\downarrow&&\,\,\downarrow j&&\,\,\downarrow h\\ G_{SC}({\mathbb A}_{F})&\stackrel{\pi}{\to}&G({\mathbb A}_{F})&\stackrel{ab_{{\mathbb A}_{F}}}{\to} &H^0_{ab}({\mathbb A}_{F};G)&\to& \prod_{v\in Val(F)} H^1(F_{v};G_{SC}) \\ \end{array}$$
Les suites horizontales sont exactes, cf. \cite{L} page 31. Soit $x\in H^0_{ab}(F;G)$, supposons $j(x)\in ab_{{\mathbb A}_{F}}(G({\mathbb A}_{F}))$. Alors l'image de $x$ dans le terme sud-est du diagramme ci-dessus est triviale. L'image de $x$ dans $H^1(F;G_{SC})$ est donc dans le noyau de $h$. Or ce noyau est trivial (th\'eor\`eme de Kneser-Harder-Chernousov, \cite{L} th\'eor\`eme 1.6.9). Donc $x$ appartient \`a l'image de $ab_{F}$ et $j(x)$ appartient \`a $ j\circ ab_{F}(G(F))$. Cela prouve (6) et ach\`eve la d\'emonstration du lemme.

\bigskip

e-mail: jean-loup.waldspurger@imj-prg.fr

 \end{document}